\documentclass[a4paper,pdftex,12pt]{article}

\usepackage{geometry}
\usepackage{times}
\usepackage{amssymb,amsmath}
\usepackage{epsfig, graphicx}

\newtheorem{Theorem}{Theorem}[section]
\newtheorem{Corollary}[Theorem]{Corollary}
\newtheorem{Lemma}[Theorem]{Lemma}
\newenvironment{Proof}
{\begin{trivlist}\item[]{{\sc Proof.}}}{\hfill{$\square$}\noindent\end{trivlist}}

\newtheorem{Definition}[Theorem]{Definition}

\title{$\mathbf{3}$-regular matchstick graphs with given girth}
\author{Sascha Kurz$^\star$ and Giuseppe Mazzuoccolo$^{\star\star}$}
\date{\small $^\star$Fakult\"at f\"ur Mathematik, Physik und Informatik, Universit\"at Bayreuth, Germany, sascha.kurz@uni-bayreuth.de\\[2mm]
      $^{\star\star}$Dipartimento di Matematica, Universitat di Modena e Reggio Emilia, via Campi 213/B, I-41100 Modena, Italy,
      mazzuoccolo.giuseppe@unimore.it}

\begin{document}

\maketitle

{\small
\noindent \textbf{Abstract:}
  We consider $3$-regular planar matchstick graphs, i.~e.{} those which have a planar embedding such that all
  edge lengths are equal, with
  given girth $g$. For girth $3$ it is known that such graphs exist if and only if the number of vertices $n$
  is an even integer larger or
  equal to $8$. Here we prove that such graphs exist for girth $g=4$ if and only if $n$ is even and at least $20$.
  We provide an example for girth $g=5$ consisting of $180$ vertices.
  }

\section{Introduction}
\noindent
In August of 1986, a special conference on recreational mathematics was held at the University of Calgary to celebrate the founding of the Strens Collection. Leading practitioners of recreational mathematics from around the world gathered in Calgary to share with each other the joy and spirit of play that is to be found in recreational mathematics, see \cite{0828.00001}. Heiko Harborth, one of the presenters, took the chance to insist that ``Matchsticks are the cheapest and simplest objects for puzzles which can be both challenging and mathematical''.

More than 20 years later our knowledge on matchstick graphs, i.~e.{} noncrossing arrangements of matchsticks, is still very limited. It seems to be hard to obtain rigid mathematical results about them.

One of those puzzles asks for the complete $r$-regular matchstick graph, see definitions \ref{def_1} and \ref{def_2}, with the minimum number of vertices, see e.~g.{} \cite{0591.52012}. Checking all possibilities, as one can do in most puzzles, is not that easy for matchstick graphs. There are some computational and mathematical obstructions. These geometric graphs may be flexible, i.~e.{} one can not determine an up to isomorphisms finite list of sets of coordinates for the vertices. Indeed the smallest complete $3$-regular matchstick graph is flexible. Also in the case where the graphs are rigid (i.~e.{} not flexible) it can be a hard task to determine the coordinates of the vertices. As an example we refer the interested reader to \cite{gerbracht} where the coordinates of the so-called Harborth graph. i.~e.{} the smallest known $4$-regular matchstick graph, were determined. Some minimal polynomials of vertex coordinates have a degree of $22$. And indeed testing whether a given planar graph can be realized in the Euclidean plane is an NP-hard problem, see \cite{pre05493581,0744.05053}.

What is known about this specific puzzle? Due to the Eulerian polyhedron formula there can not exist (finite) complete $r$-regular matchstick graphs for $r\ge 6$. For $r=1$ the only example is a single edge and for $r=2$ the only examples are circles $C_n$ for $n\ge 3$. In the next case the smallest possible number $n$ of vertices of a complete $3$-regular matchstick graph is $8$. We leave it as an easy but entertaining exercise to the reader to proof that complete $3$-regular matchstick graphs exist if and only if $n\ge 8$ is an even number and that there is at least one example. For $r=4$ we have already mentioned that the smallest known example is the so called Harborth graph consisting of $n=52$ vertices.
Very recently one of the authors proves the non-existence of a (finite) complete $5$-regular matchstick graph, see \cite{no_five_regular}, and indeed a lot of non-trivial mathematics is involved. So only the case $r=4$ remains open, but so far it seems to be out of reach. A generalization of this problem were more than one possible edge lengths is allowed, is considered in \cite{1063.05036}. In \cite{0629.52012} the author considers complete $r$-regular graphs where the edges have unit length but are allowed to cross.

In this article we consider another matchstick puzzle -- complete $3$-regular matchstick graphs with given girth and minimum number of vertices. At the end of 2005 Erich Friedman posed this problem on his ``Math magic''-homepages\footnote[1]{http://www.stetson.edu/~efriedma/mathmagic/1205.html}. He was especially interested in an example for girth $g=4$. Very soon Gavin Theobald found such an example consisting of $40$ vertices, which was beaten by an example of one of the authors consisting of only $32$ vertices in 2006.

Locating these examples is a creative and recreational task. For some matchstick problems constructing the minimal example can be quite challenging. But the really hard task is to rigidly prove that no smaller example can exist. Here we want to demonstrate that it is possible, with admittedly quite some effort, to rigidly solve a matchstick puzzle where the minimal answer has $20$ vertices. Nevertheless we aim to solve a very specific puzzle we try to present the underlying ideas and techniques from a more general point of view.

In the remaining part of this article we prove that a complete $3$-regular matchstick graph with girth $4$ exists if and only if $n$ is an even number greater or equal to $20$. We give an example of a complete $3$-regular matchstick graph with girth $5$ consisting of $180$ vertices and provide a first lower bound on the minimum number of necessary vertices. As a simple consequence from the Eulerian polyhedron formula there are no complete $3$-regular matchstick graphs with girth at least $6$ (besides from the infinite honeycomb lattice), see Lemma~\ref{lemma_A_i_sum} and Equation~(\ref{eq_A_i_sum}). Now let us go into the details.

\begin{Definition}
  \label{def_1}
  An \textbf{(incomplete) $\mathbf{r}$-regular matchstick graph} $\mathcal{M}$ consists of a graph $G=(V,E)$
  and an embedding $f:V\rightarrow\mathbb{R}^2$ in the plane which fulfill the following conditions:
  \begin{itemize}
    \item[(1)] $G$ is a connected planar graph.
    \item[(2)] The nodes on the outer face of $\mathcal{M}$ all have
               degree at most $r$ and all other nodes have degree exactly $r$.
    \item[(2)] If $\{i,j\}\in E$ then we have $\Vert f(i),f(j)\Vert_2=1$, where $\Vert x,y\Vert_2$
               denotes the Euclidean distance between the vectors $x$ and $y$.
    \item[(3)] For $i\neq j$ we have $f(i)\neq f(j)$.
    \item[(4)] If $\{i_1,j_1\},\{i_2,j_2\}\in E$ for pairwise different $i_1,i_2,j_1,j_2\in V$
               then the line segments $\overline{f(i_1)f(j_1)}$ and $\overline{f(i_2)f(j_2)}$
               do not have a common point.
  \end{itemize}
\end{Definition}

So in other words an $r$-regular matchstick graph $\mathcal{M}$ is an embedded planar graph where the inner vertices have degree $r$ and the edges are straight line segments of length $1$.

\begin{Definition}
  \label{def_2}
  We call an $r$-regular matchstick graph \textbf{complete} if all nodes on the outer face of
  $\mathcal{M}$ have degree exactly $r$.
\end{Definition}

\section{Basic definitions and parameters of matchstick graphs}
\label{sec_basic}

\noindent
In this section we introduce some parameters of matchstick graphs and prove some necessary conditions and restrictions on these parameters.

\begin{Definition}
  For the number $|V|$ of vertices of $\mathcal{M}$ we introduce the abbreviation $n(\mathcal{M})$.
  By $\mathcal{K}(\mathcal{M})$ we denote the set of vertices which is situated on the outer face of
  $\mathcal{M}$ and by $\mathcal{I}(\mathcal{M})$ we denote the set of the remaining vertices. For the
  cardinality of $\mathcal{K}(\mathcal{M})$ we introduce the notation $k(\mathcal{M})$. By
  $\tau(\mathcal{M})$ we denote the quantity
  $r\cdot k(\mathcal{M})-\sum\limits_{v\in\mathcal{K}(\mathcal{M})}\delta(v)$, where $\delta(v)$ denotes the
  degree of vertex $v$. By $A_i(\mathcal{M})$ we denote the number of faces of $\mathcal{M}$ which are
  $i$-gons. Here we also count the outer face.
\end{Definition}

Whenever it is clear from the context which matchstick graph $\mathcal{M}$ is meant we only write $n$, $\mathcal{K}$, $\mathcal{I}$, $k$, $\tau$, $A_i$ instead of $n(\mathcal{M})$, $\mathcal{K}(\mathcal{M})$, $\mathcal{I}(\mathcal{M})$, $k(\mathcal{M})$, $\tau(\mathcal{M})$, $A_i(\mathcal{M})$.

One of the basic tools for planar graphs is the Eulerian polyhedron formula which leads to the following lemma:

\begin{Lemma}
  \label{lemma_A_i_sum}
  $$
    \sum_{i=3}^\infty (2i-ri+2r)\cdot A_i=4r-2\tau.
  $$
\end{Lemma}
\begin{Proof}
  Because every edge belongs to two faces and every vertex is part of $r$ faces except some of the
  vertices of the outer face we have
  $$
    2\cdot|E|=\sum_{i=3}^\infty i\cdot A_i\quad\text{and}\quad r\cdot |V|=\sum_{i=3}^\infty i\cdot
    A_i\,+\,\tau.
  $$
  Using $|F|:=\sum\limits_{i=3}^\infty A_i$ for the number of faces and plugging these equations into the Eulerian
  polyhedron formula $|V|-|E|+|F|=2$ we obtain the stated formula.
\end{Proof}

In the remaining part of this article we will focus an the case $r=3$. Inserting $r=3$ in Lemma \ref{lemma_A_i_sum} yields
\begin{equation}
  3A_3+2A_4+A_5-A_7-2A_8-3A_9-\dots =12-2\tau\label{eq_A_i_sum}.
\end{equation}
Thus there are no (finite) complete $3$-regular matchstick graphs with girth larger than $5$.

\begin{Lemma}
  \label{lemma_V_bound}
  For a $3$-regular matchstick graph we have
  $$
    n=|V|=2|F|-4+\tau.
  $$
\end{Lemma}
\begin{Proof}
  Direct conclusion from the Eulerian polyhedron formula.
\end{Proof}
We remark that for $r=3$ we have $n\equiv \tau \pmod 2$.

In order to exclude graphs and embeddings not using the embedding function $f:V\rightarrow\mathbb{R}^2$ we may utilize some area and perimeter arguments. The outer face is a $k$-gon and their exist an upper bound on the area of a $k$-gon with side lengths $1$. Since everything else must be inside this $k$-gon we get some restrictions on the parameters of matchstick graphs.

\begin{Definition}
  For 
  $a_1,\dots a_r$ we define $A_{\text{max}}(a_1,\dots,a_r)\in
  \mathbb{R}_{\ge 0}\cup\{-1\}$ by:
  \begin{itemize}
    \item[(1)] For all polygons $P$ with side lengths $a_1,\dots,a_r$ in an arbitrary ordering we
               have that the area of $P$ is at most $A_{\text{max}}(a_1,\dots,a_r)$.
    \item[(2)] For every $\varepsilon>0$ there exists a polygon $P$ consisting of $r$ sides with
               lengths $a_1,\dots,a_r$ in a suitable ordering, where $P$ has an area of at least
               $A_{\text{max}}(a_1,\dots,a_r)-\varepsilon$.
    \item[(3)] If no polygon $P$ with side lengths $a_1,\dots,a_r$ exists then we set
               $A_{\text{max}}(a_1,\dots,a_r)=-1$.
  \end{itemize}
\end{Definition}

The determinantion of $A_{\text{max}}(a_1,\dots,a_r)$ and the characterization of the extremal examples is a well know problem in plane geometry. If all edge lengths are equal the extremal examples are the regular $r$-gons:

\begin{Lemma}
  $$
    A_{\text{max}}
    \underset{r\text{ one's}}{\underbrace{(1,\dots,1)}}=\frac{r}{4}\cdot\cot\left(\frac{\pi}{r}\right).
  $$
\end{Lemma}

For the ease of notation we use $A_{\text{max}}(r)$ instead of a sequence of $r$ ones. We remark that such an area argument is extremely useful if we can conclude the presence of several triangles, which each have an area of $\frac{\sqrt{3}}{4}$. Since in our problem we assume $A_3=0$ we also need some bounds on the minimum possible area of equilateral $i$-gons with $i\ge 4$.

\begin{table}[ht]
  \begin{tabular}{|rr|rr|rr|rr|rr|rr|}
    \hline
    k & $\frac{A_{\text{max}}(k)}{A_{\text{max}}(3)}$ & k & $\frac{A_{\text{max}}(k)}{A_{\text{max}}(3)}$ &
    k & $\frac{A_{\text{max}}(k)}{A_{\text{max}}(3)}$ & k & $\frac{A_{\text{max}}(k)}{A_{\text{max}}(3)}$ &
    k & $\frac{A_{\text{max}}(k)}{A_{\text{max}}(3)}$ \\
    \hline
     3 &   1.000 &  4 &   2.310 &  5 &   3.974 &  6 &   6.000 &  7 &   8.393 \\
     8 &  11.151 &  9 &  14.277 & 10 &  17.770 & 11 &  21.630 & 12 &  25.857 \\
    \hline
  \end{tabular}
  \caption{An upper bound for the number of equilateral triangles inside an equilateral $k$-gon.}
  \label{table_max_area}
\end{table}

\begin{Definition}
   For 
   $a_1,\dots a_r$ we define $A_{\text{min}}(a_1,\dots,a_r)\in
  \mathbb{R}_{\ge 0}\cup\{-1\}$ by:
  \begin{itemize}
    \item[(1)] For all polygons $P$ with side lengths $a_1,\dots,a_r$ in an arbitrary ordering we
               have that the area of $P$ is at least $A_{\text{min}}(a_1,\dots,a_r)$.
    \item[(2)] For every $\varepsilon>0$ there exists a polygon $P$ consisting of $r$ sides with
               lengths $a_1,\dots,a_r$ in a suitable ordering, where $P$ has an area of at most
               $A_{\text{max}}(a_1,\dots,a_r)+\varepsilon$.
    \item[(3)] If no polygon $P$ with side lengths $a_1,\dots,a_r$ exists then we set
               $A_{\text{max}}(a_1,\dots,a_r)=-1$.
  \end{itemize}
\end{Definition}

Concerning the minimum area of an equilateral $r$-gon in \cite{0968.51014} the authors have proved:
\begin{Theorem}
  \label{thm_min_area}
  $$
    A_{\text{min}}(r):=A_{\text{min}}\underset{r\text{ one's}}{\underbrace{(1,\dots,1)}}=
    \left\{\begin{array}{rcl}\frac{\sqrt{3}}{4}&:&r\equiv 1\pmod 2,\, r\ge 3,\\0&:&r\equiv 0\pmod 2,\,r\ge 4,\\-1&:&\text{else}.\end{array}\right.
  $$
\end{Theorem}

So e.~g.{} an inner pentagon in a matchstick graph has an area of at least $\frac{\sqrt{3}}{4}$. Unfortunately the area of a quadrangle may be arbitrary small. So in the case of girth $4$ we see no easy way to utilize area arguments.

We can also utilize information on the perimeter of subconfigurations to deduce the impossibility of some cases.
\begin{Lemma}
  \label{lemma_perimeter}
  If $P$ is a polygon inside another polygon $Q$, where $P\neq Q$, then the perimeter of $Q$ has to be strictly larger than
  the perimeter of the convex hull of $P$.
\end{Lemma}
So e.~g.{} we have a proof for the obvious fact that no equilateral quadrangle can contain another equilateral quadrangle of the same edge length.

\begin{Lemma}
  \label{lemma_neighbored_angles}
  If $\alpha$ and $\beta$ are two neighbored angles of an equilateral $k$-gon, then $\alpha+\beta>\frac{\pi}{2}$ holds.
\end{Lemma}
\begin{Proof}
  We assume w.l.o.g.{} $\alpha\le\beta$. For minimal $\alpha+\beta$ we could assume that the two arms of the neighbored angles touch,
  so that we have an isosceles triangle with $\alpha+2\beta>\pi$. For $\beta\le\frac{\pi}{2}$ we have
  $\alpha+\beta>\pi-\beta\ge\frac{\pi}{2}$.
\end{Proof}

\section{3-regular matchstick graphs with girth 4}

\noindent
In this section we prove our main theorem, i.~e.{} $3$-regular matchstick graph with girth $4$ exist if and only if the number of their vertices is even and greater or equal to $20$.

\begin{Lemma}
  \label{lemma_no_point_in_quadrangle}
  The only possible bridgeless connected $3$-regular matchstick graph with $k=4$ and girth $g=4$ is a simple quadrangle.
\end{Lemma}
\begin{Proof}
  Due to Lemma~\ref{lemma_perimeter} we have $A_4=1$ (the outer face) for $k=4$. If $\mathcal{M}$ contains further vertices or edges besides
  the four vertices and four edges of the outer face then we have $0\le\tau\le 3$. From Equation~(\ref{eq_A_i_sum}) we obtain
  $A_5\ge 12-2\cdot 1-2\cdot3=4$, which is a contradiction to Theorem \ref{thm_min_area} and Table~\ref{table_max_area}.
\end{Proof}

In other words Lemma \ref{lemma_no_point_in_quadrangle} says that there is no vertex and no edge inside a quadrangle of a bridgeless connected $3$-regular matchstick graph with girth at least $4$.

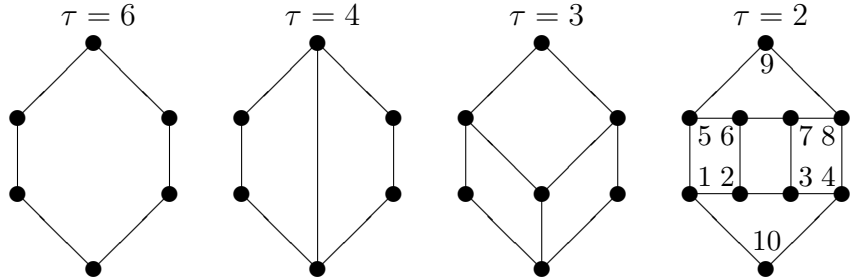
\begin{figure}[htp]
  \begin{center}
  \setlength{\unitlength}{1.0cm}
  \begin{picture}(2,3.5)
    \put(1,0){\circle*{0.2}}
    \put(0,1){\circle*{0.2}}
    \put(2,1){\circle*{0.2}}
    \put(0,2){\circle*{0.2}}
    \put(2,2){\circle*{0.2}} 
    \put(1,3){\circle*{0.2}}
    \put(1,0){\line(-1,1){1}}
    \put(1,0){\line(1,1){1}}
    \put(1,3){\line(-1,-1){1}}
    \put(1,3){\line(1,-1){1}}
    \put(0,1){\line(0,1){1}}
    \put(2,1){\line(0,1){1}}
    \put(0.57,3.23){$\tau=6$}
  \end{picture}
  \quad\quad
  \begin{picture}(2,3.5)
    \put(1,0){\circle*{0.2}}
    \put(0,1){\circle*{0.2}}
    \put(2,1){\circle*{0.2}}
    \put(0,2){\circle*{0.2}}
    \put(2,2){\circle*{0.2}} 
    \put(1,3){\circle*{0.2}}
    \put(1,0){\line(-1,1){1}}
    \put(1,0){\line(1,1){1}}
    \put(1,3){\line(-1,-1){1}}
    \put(1,3){\line(1,-1){1}}
    \put(0,1){\line(0,1){1}}
    \put(2,1){\line(0,1){1}}
    \put(1,0){\line(0,1){3}}
    \put(0.57,3.23){$\tau=4$}
  \end{picture}
  \quad\quad
  \begin{picture}(2,3.5)
    \put(1,0){\circle*{0.2}}
    \put(0,1){\circle*{0.2}}
    \put(2,1){\circle*{0.2}}
    \put(0,2){\circle*{0.2}}
    \put(2,2){\circle*{0.2}} 
    \put(1,3){\circle*{0.2}}
    \put(1,1){\circle*{0.2}}
    \put(1,0){\line(-1,1){1}}
    \put(1,0){\line(1,1){1}}
    \put(1,3){\line(-1,-1){1}}
    \put(1,3){\line(1,-1){1}}
    \put(0,1){\line(0,1){1}}
    \put(2,1){\line(0,1){1}}
    \put(1,0){\line(0,1){1}}
    \put(1,1){\line(-1,1){1}}
    \put(1,1){\line(1,1){1}}
    \put(0.57,3.23){$\tau=3$}
  \end{picture}
  \quad\quad
  \begin{picture}(2,3.5)
    \put(1,0){\circle*{0.2}}
    \put(0,1){\circle*{0.2}}
    \put(2,1){\circle*{0.2}}
    \put(0,2){\circle*{0.2}}
    \put(2,2){\circle*{0.2}} 
    \put(1,3){\circle*{0.2}}
    \put(0.66,1){\circle*{0.2}}
    \put(1.33,1){\circle*{0.2}}
    \put(0.66,2){\circle*{0.2}}
    \put(1.33,2){\circle*{0.2}}
    \put(1,0){\line(-1,1){1}}
    \put(1,0){\line(1,1){1}}
    \put(1,3){\line(-1,-1){1}}
    \put(1,3){\line(1,-1){1}}
    \put(0,1){\line(0,1){1}}
    \put(2,1){\line(0,1){1}}
    \put(0.66,1){\line(0,1){1}}
    \put(1.33,1){\line(0,1){1}}
    \put(0,1){\line(1,0){2}}
    \put(0,2){\line(1,0){2}}
    \put(0.1,1.1){\small{$1$}}
    \put(0.4,1.1){\small{$2$}}
    \put(1.43,1.1){\small{$3$}}
    \put(1.73,1.1){\small{$4$}}
    \put(0.1,1.65){\small{$5$}}
    \put(0.4,1.65){\small{$6$}}
    \put(1.43,1.65){\small{$7$}}
    \put(1.73,1.65){\small{$8$}}
    \put(0.92,2.6){\small{$9$}}
    \put(0.82,0.25){\small{${10}$}}
    \put(0.57,3.23){$\tau=2$}
  \end{picture}
  \caption{Planar matchstick graphs with $k=6$.}
  \label{fig_matchstick_graphs_1}
  \end{center}
\end{figure}

\begin{Lemma}
  \label{lemma_no_inner_4_4_4}
  If a bridgeless connected $3$-regular matchstick graph $\mathcal{M}$ with girth at least $4$ contains an inner point $v$ whose
  three adjacent faces all are quadrangles, then $k\ge 6$ and $k=6$ is possible only if $v$ is the unique inner point (see the third
  graph of Figure \ref{fig_matchstick_graphs_1}).
\end{Lemma}
\begin{Proof}
  Let $\alpha$, $\beta$, and $\gamma$ be the three angles of the quadrangles at $v$. In an equilateral quadrangle we have
  $\alpha,\beta,\gamma\in(0,\pi)$. Since $\alpha+\beta+\gamma=2\pi$ we have $\min(\alpha+\beta,\alpha+\gamma,\beta+\gamma)>\pi$.
  Now we consider the matchstick graph $\mathcal{M}'$ consisting of $v$ and the three adjacent quadrangles. The outer face of $\mathcal{M}'$
  is a hexagon with inner angles $\alpha,\beta,\gamma,2\pi-\alpha-\beta,2\pi-\alpha-\gamma,2\pi-\beta-\gamma$. Thus all angles are smaller than
  $\pi$ and the outer face is a convex hexagon with edges of unit length. Now the statement follows from Lemma~\ref{lemma_perimeter}.
\end{Proof}

\begin{figure}[htp]
\centering
\includegraphics[width=8cm]{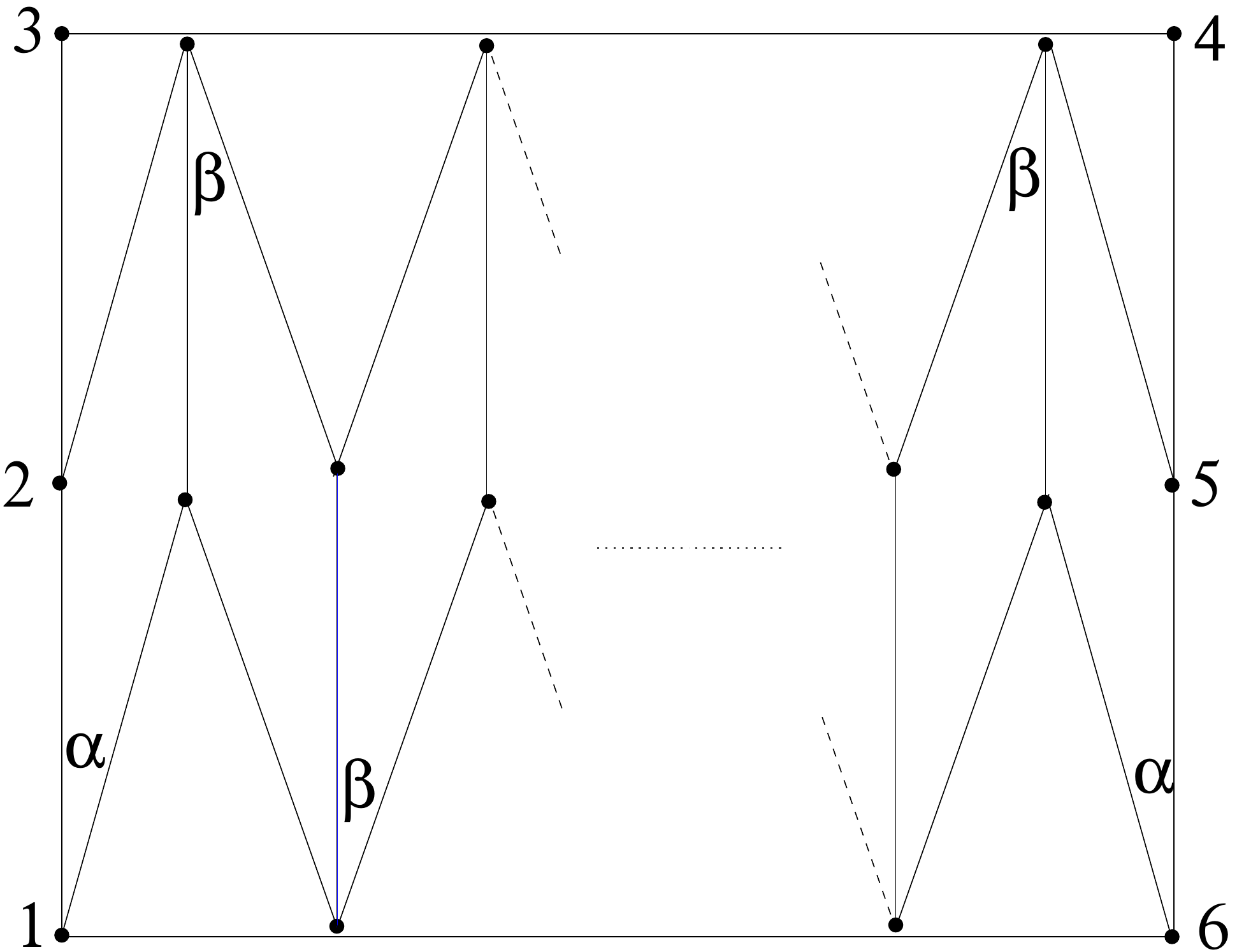}
\caption{Arbitrarily large matchstick graphs for $k=6$.}
 \label{fig_k_6}
\end{figure}

We would like to remark that we have tried some attempts in order to prove that the simple pentagon is the only possibility for $k=5$, girth $4$, and arbitrary $\tau$, but every time we have ended up in confusing case differentiations.
On the other hand there exist arbitrarily large $3$-regular matchstick graphs with $k=6$,
a possible construction is given in Figure \ref{fig_k_6}: one can make a suitable choice of angles $\alpha$ and $\beta$ in such a way that $0<\alpha<\beta<\frac{\pi}{2}$, the distance between vertices $2$ and $5$ is equal to $1$, and the vertices $1,\dots,6$ form a $2\times 1$ rectangle. 
To be more precisely, if the coordinates of vertex $1$ are $\begin{pmatrix}0 & 0\end{pmatrix}$, the coordinates of vertex $2$ are $\begin{pmatrix}0 & 1\end{pmatrix}$, and we have a chain of $2n$ quadrangles, then coordinates of vertices $5$ and $6$ are given by $\begin{pmatrix}2\sin\alpha+(2n-2)\sin\beta & 1 \end{pmatrix}$ and $\begin{pmatrix}2\sin\alpha+(2n-2)\sin\beta & 0 \end{pmatrix}$, respectively. If $\alpha<\beta$ the $y$-coordinates of all vertices except vertex $1,\ldots,6$ are contained in the open interval $(0,2)$. Choosing a pair $\alpha,\beta$, such that the $x$-coordinate of vertices $4,5,6$ is $1$, gives the construction.

\begin{figure}[htp]
\begin{center}
  \setlength{\unitlength}{3.0cm}
  \begin{picture}(0.8,1.3)
    \put(0.3007057995 , -0.9537169507){\circle*{0.025}} 
    \put(0.3007057995 , 0.0462830493){\circle*{0.025}}  
    \put(0.8742822351 , -0.7728689953){\circle*{0.025}} 
    \put(0.8742822351 , 0.2271310053){\circle*{0.025}}  
    \put(0.0000000000 , 0.0000000000){\circle*{0.025}}  
    \put(0.0000000000 , 1.0000000000){\circle*{0.025}}  
    \put(0.5735764363 , 0.1808479557){\circle*{0.025}}  
    \put(0.5735764363 , 1.1808479560){\circle*{0.025}}  
    \put(-0.3918153149 , 0.9200438882){\circle*{0.025}} 
    \put(1.266097554 , -0.6929128843){\circle*{0.025}}  
    \put(-0.4618153149 , 0.8800438882){\small{9}}
    \put(1.296097554 , -0.7329128843){\small{10}}
    \qbezier(0.3007057995 , -0.9537169507)(0.3007057995 , -0.9537169507)(0.0000000000 , 0.0000000000) 
    \qbezier(0.3007057995 , 0.0462830493)(0.3007057995 , 0.0462830493)(0.0000000000 , 1.0000000000)   
    \qbezier(0.8742822351 , -0.7728689953)(0.8742822351 , -0.7728689953)(0.5735764363 , 0.1808479557) 
    \qbezier(0.8742822351 , 0.2271310053)(0.8742822351 , 0.2271310053)(0.5735764363 , 1.1808479560)   
    \qbezier(0.0000000000 , 0.0000000000)(0.0000000000 , 0.0000000000)(0.0000000000 , 1.0000000000)   
    \qbezier(0.0000000000 , 1.0000000000)(0.0000000000 , 1.0000000000)(0.5735764363 , 0.1808479557)   
    \qbezier(0.5735764363 , 0.1808479557)(0.5735764363 , 0.1808479557)(0.5735764363 , 1.1808479560)   
    \qbezier(0.3007057995 , -0.9537169507)(0.3007057995 , -0.9537169507)(0.3007057995 , 0.0462830493) 
    \qbezier(0.3007057995 , 0.0462830493)(0.3007057995 , 0.0462830493)(0.8742822351 , -0.7728689953)  
    \qbezier(0.8742822351 , -0.7728689953)(0.8742822351 , -0.7728689953)(0.8742822351 , 0.2271310053) 
    \qbezier(0.3007057995 , -0.9537169507)(0.3007057995 , -0.9537169507)(1.266097554 , -0.6929128843) 
    \qbezier(0.8742822351 , 0.2271310053)(0.8742822351 , 0.2271310053)(1.266097554 , -0.6929128843)   
    \qbezier(0.0000000000 , 0.0000000000)(0.0000000000 , 0.0000000000)(-0.3918153149 , 0.9200438882)  
    \qbezier(0.5735764363 , 1.1808479560)(0.5735764363 , 1.1808479560)(-0.3918153149 , 0.9200438882)  
  \end{picture}\\[24.5mm]
  \caption{A matchstick graph for $k=6$ and $\tau=2$ consisting of $10$ vertices.}
  \label{fig_k_6_tau_2}
\end{center}
\end{figure}
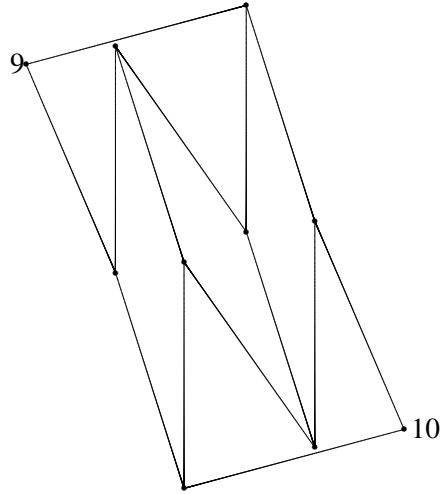

\begin{Lemma}
  The fourth (embedded) planar graphs of Figure \ref{fig_matchstick_graphs_1} can be realized with edges of unit length, so that
  the outer face is not convex.
\end{Lemma}
\begin{Proof}
  Let $\alpha$ be the outer angle of the upper pentagon between vertices $5$, $6$, and $7$. By $\beta$ we denote the neighbored
  outer angle between vertices $6$, $7$, and $8$.

  We remark that there is a small range for the possible values of $\alpha$ and $\beta$, but e.~g.{} $\alpha=\frac{35}{180}\cdot\pi$,
  $\beta=2\pi-\alpha=\frac{325}{180}\cdot\pi$ leads to a valid matchstick configuration. We set $p_5 = \begin{pmatrix}0&0\end{pmatrix}$,
  $p_6 = \begin{pmatrix}0&1\end{pmatrix}$, and let the quadrangles bisect $\alpha$ and $\beta$:
  \begin{eqnarray*}
    p_1 &=& \begin{pmatrix}\sin\left(\frac{7\pi}{72}\right)&-\cos\left(\frac{7\pi}{72}\right)\end{pmatrix}\\
        &\approx& \begin{pmatrix}0.3007057995 &  -0.9537169507\end{pmatrix}\\
    p_2 &=& \begin{pmatrix}\sin\left(\frac{7\pi}{72}\right)&1-\cos\left(\frac{7\pi}{72}\right)\end{pmatrix}\\
        &\approx&\begin{pmatrix}0.3007057995 &  0.0462830493\end{pmatrix}\\
    p_3 &=& \begin{pmatrix}\sin\left(\frac{7\pi}{36}\right)-\sin\left(\frac{79\pi}{72}\right)&
           1-\cos\left(\frac{7\pi}{36}\right)+\cos\left(\frac{79\pi}{72}\right)\end{pmatrix}\\
           &\approx&\begin{pmatrix}0.8742822351 & -0.7728689953\end{pmatrix}\\
    p_4 &=& \begin{pmatrix}\sin\left(\frac{7\pi}{36}\right)-\sin\left(\frac{79\pi}{72}\right)&
           2-\cos\left(\frac{7\pi}{36}\right)+\cos\left(\frac{79\pi}{72}\right)\end{pmatrix}\\
           &\approx&\begin{pmatrix}0.8742822351 & 0.2271310053\end{pmatrix}\\
    p_5 &=& \begin{pmatrix}0&0\end{pmatrix}\\&\approx& \begin{pmatrix}0.0000000000&0.0000000000\end{pmatrix}\\
    p_6 &=& \begin{pmatrix}0&1\end{pmatrix}\\&\approx& \begin{pmatrix}0.0000000000&1.0000000000\end{pmatrix}\\
    p_7 &=& \begin{pmatrix}\sin\left(\frac{7\pi}{36}\right)&1-\cos\left(\frac{7\pi}{36}\right)\end{pmatrix}\\
           &\approx&\begin{pmatrix}0.5735764363&0.1808479557\end{pmatrix}\\
    p_8 &=& \begin{pmatrix}\sin\left(\frac{7\pi}{36}\right)&
           2-\cos\left(\frac{7\pi}{36}\right)\end{pmatrix}\\&\approx&
           \begin{pmatrix}0.5735764363&1.1808479560\end{pmatrix}\\
    p_9&\approx&\begin{pmatrix}-0.3918153149&0.9200438882\end{pmatrix}\\
    p_{10}&\approx&\begin{pmatrix}1.266097554& -0.6929128843 \end{pmatrix}
  \end{eqnarray*}
  The coordinates of $p_9$ and $p_{10}$ are solutions of a quadratic equation, whose analytical form can simply be
  obtained by using an arbitrary computer algebra package. A well-proportioned drawing of this example can be found
  in Figure \ref{fig_k_6_tau_2}.

  It can be easily checked, that the edges do not cross and that the outer hexagon is not convex.
\end{Proof}

\begin{figure}[htp]
\begin{center}
  \includegraphics[width=5cm]{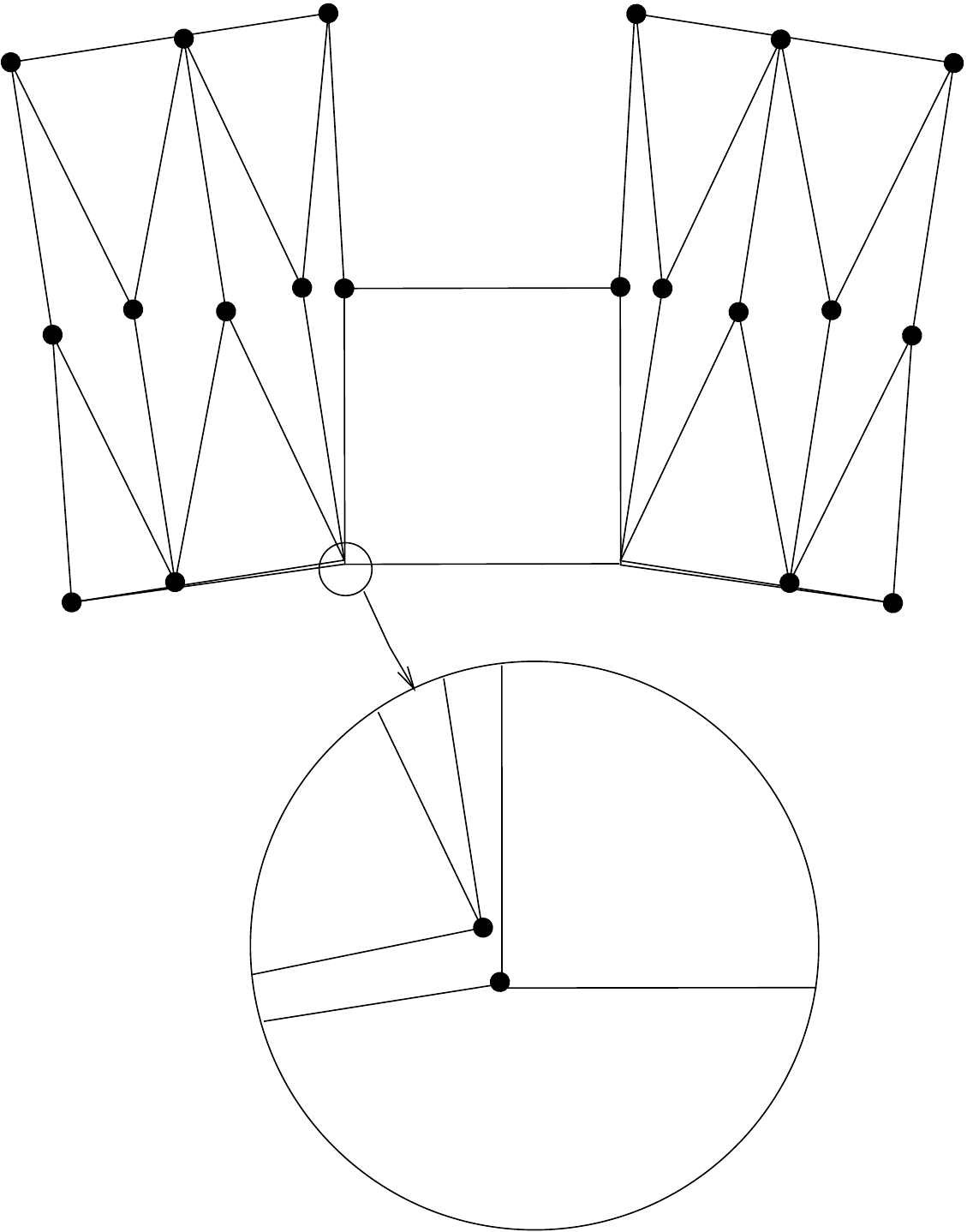}
  \caption{Construction for complete $3$-regular matchstick graphs with girth $4$ for $n\ge 22$.}
  \label{fig_constr_22_plus}
\end{center}
\end{figure}

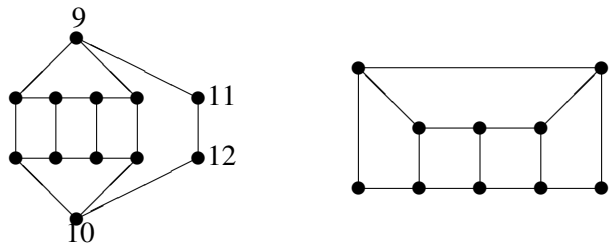
\begin{figure}[htp]
  \begin{center}
  \setlength{\unitlength}{0.8cm}
  \begin{picture}(3,3.4)
    \put(1,0){\circle*{0.2}}
    \put(0,1){\circle*{0.2}}
    \put(2,1){\circle*{0.2}}
    \put(0,2){\circle*{0.2}}
    \put(2,2){\circle*{0.2}} 
    \put(1,3){\circle*{0.2}}
    \put(0.92,3.18){\small{9}}
    \put(0.82,-0.39){\small{10}}
    \put(3.15,0.89){\small{12}}
    \put(3.15,1.89){\small{11}}
    \put(0.66,1){\circle*{0.2}}
    \put(1.33,1){\circle*{0.2}}
    \put(0.66,2){\circle*{0.2}}
    \put(1.33,2){\circle*{0.2}}
    \put(3,1){\circle*{0.2}}
    \put(3,2){\circle*{0.2}}
    \put(1,0){\line(-1,1){1}}
    \put(1,0){\line(1,1){1}}
    \put(1,3){\line(-1,-1){1}}
    \put(1,3){\line(1,-1){1}}
    \put(0,1){\line(0,1){1}}
    \put(2,1){\line(0,1){1}}
    \put(0.66,1){\line(0,1){1}}
    \put(1.33,1){\line(0,1){1}}
    \put(0,1){\line(1,0){2}}
    \put(0,2){\line(1,0){2}}
    \put(3,1){\line(0,1){1}}
    \put(3,1){\line(-2,-1){2}}
    \put(3,2){\line(-2,1){2}}
  \end{picture}\quad\quad\quad\quad\quad
  \begin{picture}(4,3)
    \put(0,0.5){\circle*{0.2}}
    \put(1,0.5){\circle*{0.2}}
    \put(2,0.5){\circle*{0.2}}
    \put(3,0.5){\circle*{0.2}}
    \put(4,0.5){\circle*{0.2}}
    \put(0,2.5){\circle*{0.2}}
    \put(1,1.5){\circle*{0.2}}
    \put(2,1.5){\circle*{0.2}}
    \put(3,1.5){\circle*{0.2}}
    \put(4,2.5){\circle*{0.2}}
    \put(0,0.5){\line(1,0){4}}
    \put(0,0.5){\line(0,1){2}}
    \put(1,0.5){\line(0,1){1}}
    \put(2,0.5){\line(0,1){1}}
    \put(3,0.5){\line(0,1){1}}
    \put(4,0.5){\line(0,1){2}}
    \put(1,1.5){\line(1,0){2}}
    \put(1,1.5){\line(-1,1){1}}
    \put(3,1.5){\line(1,1){1}}
    \put(0,2.5){\line(1,0){4}}
  \end{picture}\\[0mm]
  \caption{More planar matchstick graphs with $\tau=2$.}
  \label{fig_matchstick_2}
  \end{center}
\end{figure}

We can use the example from Figure \ref{fig_k_6_tau_2} to construct complete $3$-regular matchstick graphs with girth for $4$ consisting
of $n$ vertices for each even $n\ge 22$. Since it is not convex, we can add a path $[9,11,12,10]$ such that the points $9$, $11$, $12$, and
$10$ are in convex position, see Figure \ref{fig_matchstick_2}. If we mirror this construction, we obtain a complete $3$-regular matchstick graphs with girth for $4$ consisting of $n=22$ vertices. Replacing the edge $\{11,12\}$ by a chain of $k$ quadrangles we obtain a complete $3$-regular matchstick graphs with girth for $4$ consisting of $n=22+2k$ vertices, see Figure \ref{fig_constr_22_plus} for $k=1$.

\begin{figure}[htp]
\begin{center}
  \includegraphics[width=5cm]{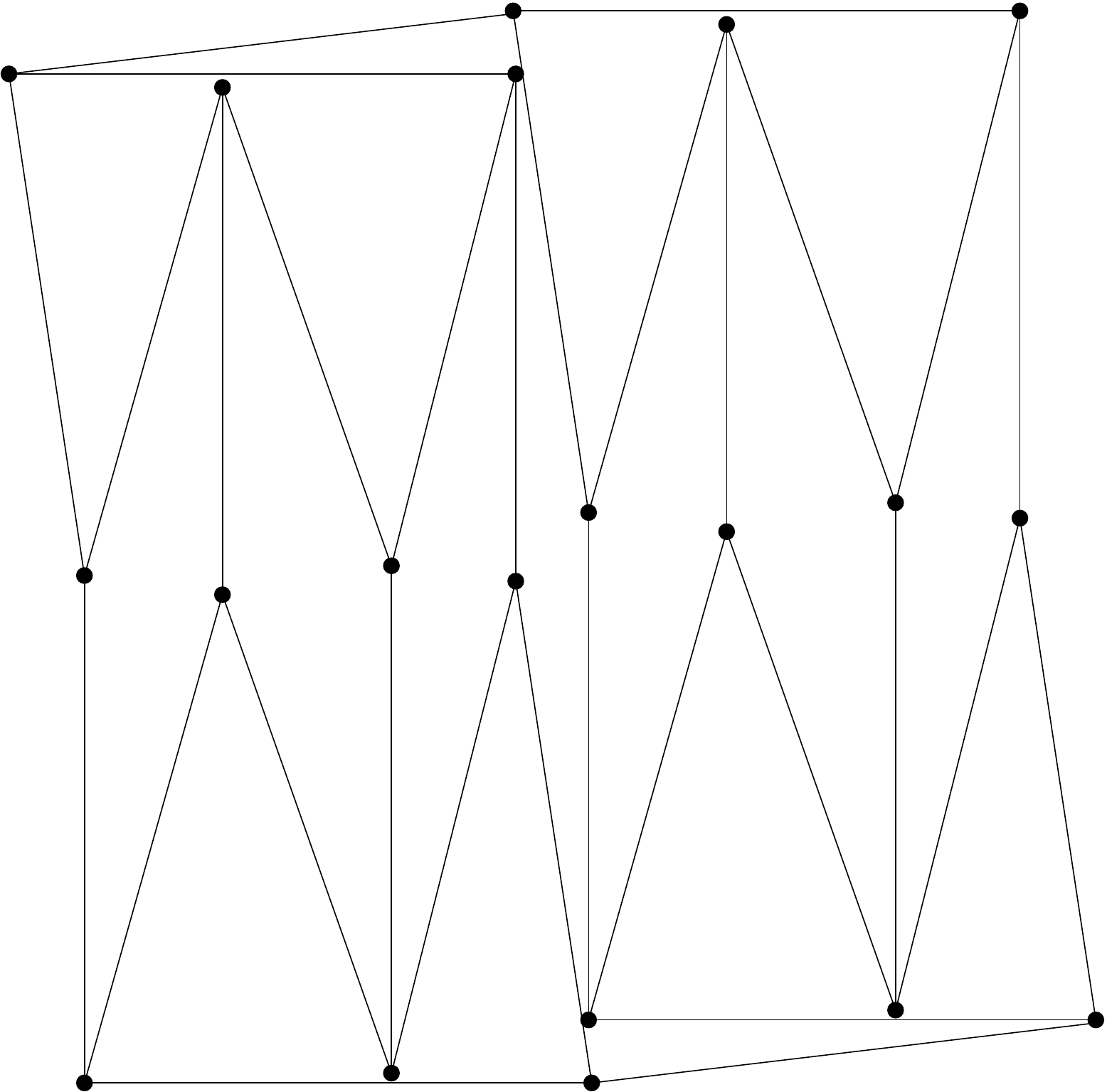}
  \caption{A complete $3$-regular matchstick graphs with girth $4$ consisting of $20$ vertices.}
  \label{fig_constr_20}
\end{center}
\end{figure}

By carefully joining two copies of the example in Figure \ref{fig_k_6_tau_2} using two additional edges, we obtain a complete $3$-regular matchstick graphs with girth $4$ consisting of $20$ vertices, see Figure \ref{fig_constr_20}.

In the remaining part of this section we prove, that there is no complete $3$-regular matstick graph with girth $4$ and fewer than $20$ vertices.

\begin{Lemma}
  \label{lemma_small}
  If $\mathcal{M}$ is a $3$-regular matchstick graph with $n\le 10$ and $1\le\tau\le 2$ then
  $\mathcal{M}$ equals either the fourth graph of Figure \ref{fig_matchstick_graphs_1} or the second graph of Figure \ref{fig_matchstick_2}.
\end{Lemma}
\begin{Proof}
  We have utilized the computer program \texttt{plantri}, see e.~g.{} \cite{pre05382101,plantri}, in order to exhaustively generate embeddings of
  planar graphs consisting of at most $10$ vertices. Checking the girth, vertex degree, $\tau$, and removing all configurations
  where Lemma \ref{lemma_no_point_in_quadrangle} or Lemma \ref{lemma_no_inner_4_4_4} can be applied, leaves only the two mentioned graphs.
\end{Proof}

In both graphs of Lemma \ref{lemma_small} we have $\tau=2$ and $n=10$. We would like to remark that a vertex-minimal $3$-regular matchstick graph with girth $4$ and $\tau=0$ is obviously $1$-connected

\begin{Theorem}
  \label{thm_connectivity}
  The connectivity and the edge-connectivity are equal in every cubic graphs.
\end{Theorem}
\begin{Proof}
  See e.~g.{} \cite[p.55]{0182.57702}.
\end{Proof}

So indeed we have that every connected cubic graph, which is not $3$-connected, decomposes into at least two connected components
after deleting at most two edges.

\begin{Lemma}
  If $\mathcal{M}$ is a complete $3$-regular matchstick graph with girth $4$, which is not $3$-connected, then we have $n\ge 20$.
\end{Lemma}
\begin{Proof}
  Let us at first assume that $\mathcal{M}$ is $1$-connected but not $2$-connected.
  Since the maximum degree is at most three every cut vertex is adjacent to a bridge. If we remove a bridge of $\mathcal{M}$
  we end up with two connected components $\mathcal{C}_1$ and $\mathcal{C}_2$ with $\tau=1$. In Lemma \ref{lemma_small} we
  have shown $n(\mathcal{C}_1),n(\mathcal{C}_2)>10$.

  If $\mathcal{M}$ is $2$-connected but not $3$-connected, then there exist two edges whose removal yields two connected
  components $\mathcal{C}_1$ and $\mathcal{C}_2$ with $\tau=2$. From Lemma \ref{lemma_small} we can conclude $n\left(\mathcal{M}\right)=
  n\left(\mathcal{C}_1\right)+n\left(\mathcal{C}_2\right)\ge 2\cdot 10=20$.
\end{Proof}

\begin{Lemma}
  \label{lemma_few_quadrangles_at_the_outer_face}
  If $\mathcal{M}$ is a $3$-connected complete $3$-regular matchstick graph with girth $4$ then the outer face is adjacent
  to at most $k(\mathcal{M})-5$ quadrangles.
\end{Lemma}
\begin{Proof}
  At first we remark that every inner quadrangle has at most one edge in common with the outer face. Now we consider the sum of the
  inner angles of the outer $k$-gon, which is $(k-2)\pi$. By $x$ we denote the number of edges of the outer face which are adjacent to
  an inner quadrangle. Due to fact that the sum of two neighbored angles in an equilateral quadrangle is $\pi$ and due to Lemma
  \ref{lemma_neighbored_angles} we obtain
  $$
    x\cdot\pi+(k-x)\cdot\frac{\pi}{2}<(k-2)\cdot\pi\quad\Leftrightarrow\quad x<k-4
  $$
  for $x<k$.
\end{Proof}

\begin{Corollary}\label{corollary_few_quadrangles}
If $\mathcal{M}$ is a $3$-connected complete $3$-regular matchstick graph with girth $4$ then $A_4\le |F|-6$.
\end{Corollary}
\begin{Proof}
Since $\mathcal{M}$ is $3$-edge-connected and $3$-regular every inner face has at most one edge in common with the outer face.
Thus due to Lemma \ref{lemma_no_point_in_quadrangle} and Lemma \ref{lemma_few_quadrangles_at_the_outer_face} there exists at least $6$ no quadrangle faces, the outer one and five of the faces adjacent to it.
\end{Proof}

\begin{Lemma}
If $\mathcal{M}$ is a $3$-connected complete $3$-regular matchstick graph with girth $4$ then we have $n\ge 16$.
\end{Lemma}
\begin{Proof}
>From Equation~(\ref{eq_A_i_sum}) we obtain $2A_4+A_5\ge 12$. Obviously $A_4+A_5\le|F|$ holds, then $A_4\ge12 -|F|$.
Due to Corollary \ref{corollary_few_quadrangles} we obtain
$$|F|-6 \ge A_4 \ge 12-|F|$$
then $|F|\ge 9$ (that is $n\ge 14$).
In the case $|F|=9$ we have $A_4 \le 3$ by Corollary \ref{corollary_few_quadrangles}.
If $A_4=3$ then from Equation~(\ref{eq_A_i_sum}) $A_5=6$, due to area arguments (see Table \ref{table_max_area}) this is not possible.
If $A_4\le 2$ then $A_4+A_5\ge12-A_4\ge10$, this is a contradiction by $A_4+A_5\le9$.
\end{Proof}

So it remains to check all $3$-connected $3$-regular planar graphs with girth $4$ consisting of $n=16$ or $18$ vertices. The exhaustive generation was again done by using the computer program \texttt{plantri}: there exists $46$ such graphs but only $23$ of them satisfy Corollary \ref{corollary_few_quadrangles} ($5$ for $n=16$ and $18$ for $n=18$). 

%

%

Lemma \ref{lemma_no_point_in_quadrangle}, Lemma \ref{lemma_no_inner_4_4_4}, Lemma \ref{lemma_few_quadrangles_at_the_outer_face}, and  Theorem \ref{thm_min_area} in combination with Table \ref{table_max_area} fortunately are sufficient to exclude all these cases so that we conclude:

\begin{Theorem}
  If $\mathcal{M}$ is a complete $3$-regular matchstick graph with girth $4$, then it contains at least $20$ vertices.
\end{Theorem}

\section{3-regular matchstick graphs with girth 5}\label{sec_girth_5}
\begin{figure}[htp]
\begin{center}
  \includegraphics[width=8cm]{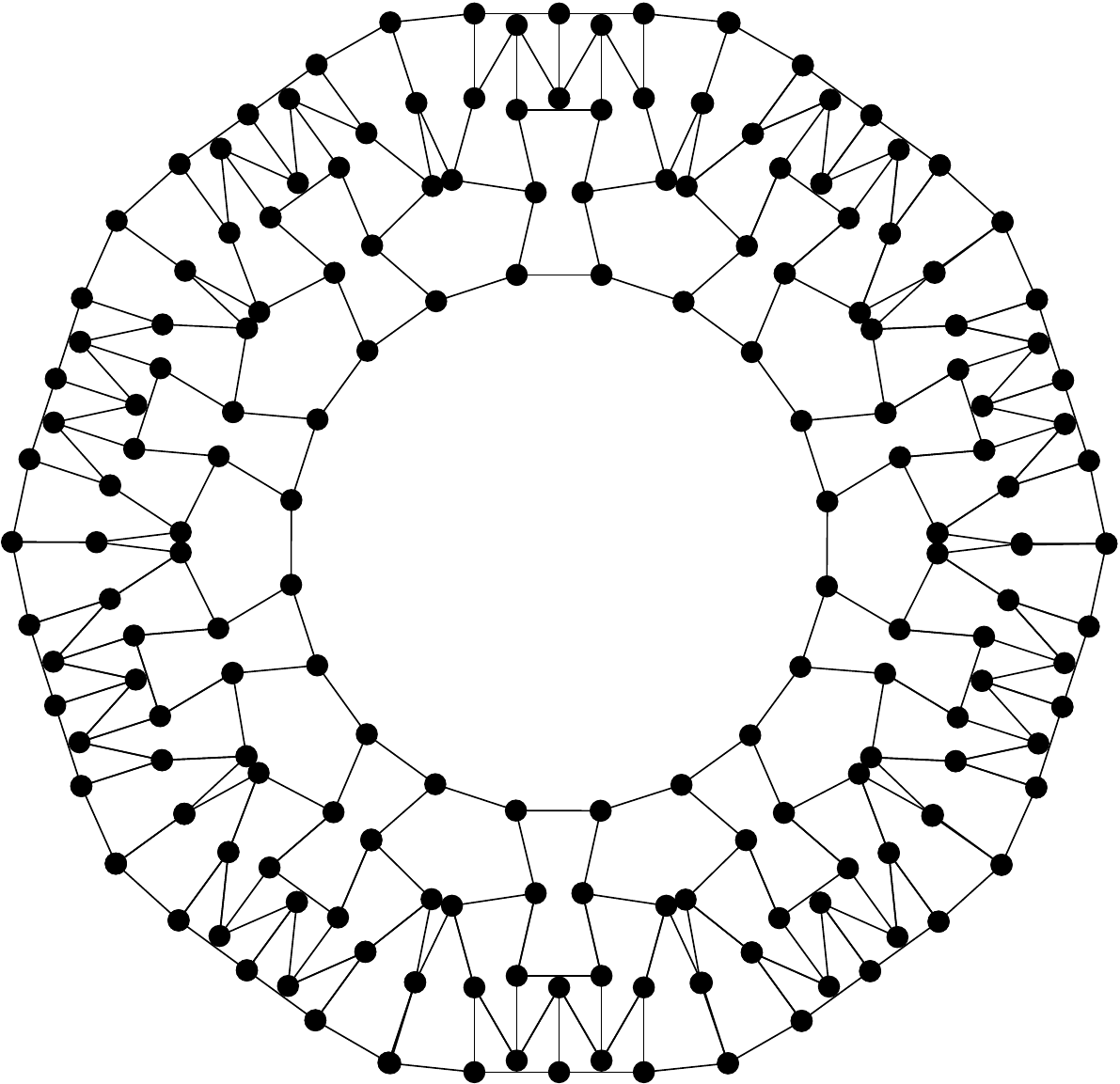}
  \caption{A complete $3$-regular matchstick graphs with girth $5$ consisting of $180$ vertices.}
  \label{fig_girth_5}
\end{center}
\end{figure}

\noindent
For complete $3$-regular matchstick graphs with girth $5$ no non-existence result, like the one from \cite{no_five_regular}, is possible, since in Figure \ref{fig_girth_5} we give an example consisting of $180$ vertices.

Applying area arguments one can obtain a first lower bound on the number of vertices of a complete $3$-regular matchstick graph with girth $5$. Inserting $\tau=0$ and $A_3=A_4=0$ into Equation~(\ref{eq_A_i_sum}) yields $A_5\ge 12$. Thus from Theorem \ref{thm_min_area}, Equation~(\ref{eq_A_i_sum}), and Table~\ref{table_max_area} we can conclude $k\ge 9$ and $A_5\ge 15$. Applying our argument again we obtain $k\ge 10$ and $A_5\ge 16$. Due to Lemma~\ref{lemma_V_bound} we have $n\ge 30$. Obviously one may easily improve this bound, but we think that the determination of the minimum example with girth $5$ may be a hard task.

\section{Conclusion and outlook}

\noindent
In Section \ref{sec_basic} we have introduced some parameters and techniques in order to provide some easy and computationally cheap certificates for proving that a given planar graph is not realizable with unit edge lengths. These methods also work for non-rigid graphs.

We have to admit that the given criteria are very far from being sufficient in general, but at least they were sufficient to completely solve a non-trivial matchstick puzzle. We would like to remark that we stumbled over the example of Figure \ref{fig_k_6_tau_2} along the way to prove the minimality of the 2006 example of Giuseppe Mazzuoccolo consisting of $32$ vertices using an exhaustive search. Indeed our first try was to prove that the configurations from Figure \ref{fig_k_6_tau_2} is not a matchstick graph. We believe that an example like those of Figure \ref{fig_constr_20} would not have been discovered by playing and puzzling with matchsticks itself. Who would try such a configuration? You would need very precisely sized and thin matchsticks.

In our opinion a lot of more needs to be done in order to provide a solid grounding for exhaustive search methods for planar geometric graphs with given side lengths. Already the case were all edges have an equal length seems to be quite hard.

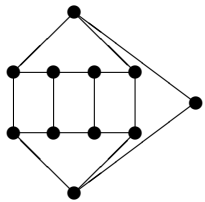
\begin{figure}[htp]
  \begin{center}
  \setlength{\unitlength}{0.8cm}
  \begin{picture}(3,3)
    \put(1,0){\circle*{0.2}}
    \put(0,1){\circle*{0.2}}
    \put(2,1){\circle*{0.2}}
    \put(0,2){\circle*{0.2}}
    \put(2,2){\circle*{0.2}} 
    \put(1,3){\circle*{0.2}}
    \put(0.66,1){\circle*{0.2}}
    \put(1.33,1){\circle*{0.2}}
    \put(0.66,2){\circle*{0.2}}
    \put(1.33,2){\circle*{0.2}}
    \put(3,1.5){\circle*{0.2}}
    \put(1,0){\line(-1,1){1}}
    \put(1,0){\line(1,1){1}}
    \put(1,3){\line(-1,-1){1}}
    \put(1,3){\line(1,-1){1}}
    \put(0,1){\line(0,1){1}}
    \put(2,1){\line(0,1){1}}
    \put(0.66,1){\line(0,1){1}}
    \put(1.33,1){\line(0,1){1}}
    \put(0,1){\line(1,0){2}}
    \put(0,2){\line(1,0){2}}
    \put(3,1.5){\line(-4,3){2}}
    \put(3,1.5){\line(-4,-3){2}}
  \end{picture}
  \caption{A planar graph which is not a matchstick graph.}
  \label{fig_forbidden_2}
  \end{center}
\end{figure}

To stimulate some research in this direction we ask the interested reader for an elegant proof or an algorithm which is relatively fast in practice to show that the planar graph from Figure \ref{fig_forbidden_2} is not a matchstick graph. Actually we do not know a general algorithm which can decide whether a given planar graph is a matchstick graph.

\end{document}